# Non-Tiles and 'Walls'- a Variant on the Heesch Problem

*Erich Friedman*
*Dept. of Mathematics, Stetson University, DeLand, FL 32720*
*erich.friedman@stetson.edu*

*R. Nandakumar*
*Amrita School of Arts and Sciences, Idapalli North, Kochi 682024, India*
*nandacumar@gmail.com*

## Introduction

The Heesch number of a 2D shape is defined as the maximum number of layers of copies of the same shape that can surround it ([1], [2], [3], [4]). A shape that tiles the plane has Heesch number infinity and the Heesch number of a shape that fails to tile the plane (a 'non-tile') gives a measure of how much it can progress towards tiling the plane. In this paper, we introduce another scheme for ranking non-tiles - using the concept of a *'wall'* and a number called *'wall thickness'* (or simply *'thickness'*).

**Definitions:** A *wall* is a simply connected region of the plane formed by infinitely many copies of a region R and dividing the plane into exactly 2 simply connected regions *A* and *B* that are a positive distance apart (intuitively, regions *A* and *B* are clearly separated by the wall and the wall has no 'cavities' inside it).

A wall *W* has *thickness* 2 (generally *n*) if it is the union of two 2 (*n*) walls which only share a boundary but is not the union of 3(*n+1*) walls. Obviously, a region that can form a wall of infinite thickness tiles the plane. For every non-tile R, we associate a *thickness number* – the thickness of the thickest wall that can be formed with copies of R.

## Examples – Walls and Thickness Numbers

1. All polygonal regions can form walls of thickness 1. But most non-tiling regions fail to form even walls of thickness 2. Figure 1 shows a very simple example of a non-tile - a circular disc with its boundary 'chopped' at the two opposite sides - that cannot form walls of thickness greater than 1, so its thickness number is 1.

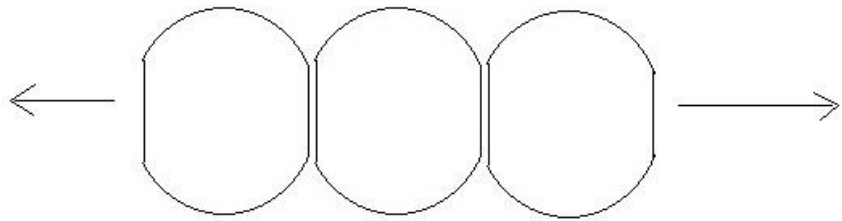

*Figure 1*

2. A non-tile which can form a wall of thickness 2 (but not more) is the convex region formed by attaching a square to a semicircle of diameter equal to the side of the square (figure 2).

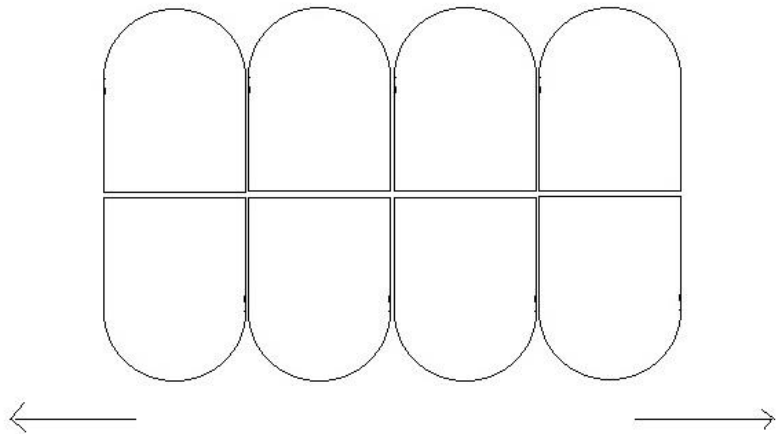

*Figure 2*

3. The Heesch Pentagon, with Heesch number 1, can form a wall of thickness 2 by repeating the pattern shown in figure 3 in the horizontal direction. Note that this arrangement is very different from the arrangement of copies of this shape that shows Heesch number 1 (see [4]).

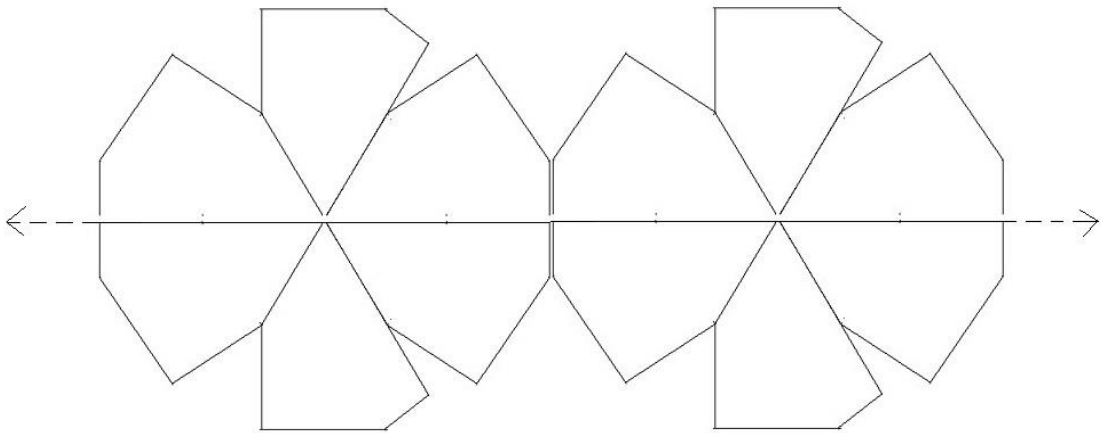

*Figure 3*

4. Figure 4 shows a non-convex region - a 5x7 rectangle with three 1x1 additions and one 1x1 hole – with Heesch number 1([2]). But it can be easily seen that this region can form walls of only thickness 1- its thickness number is 1.

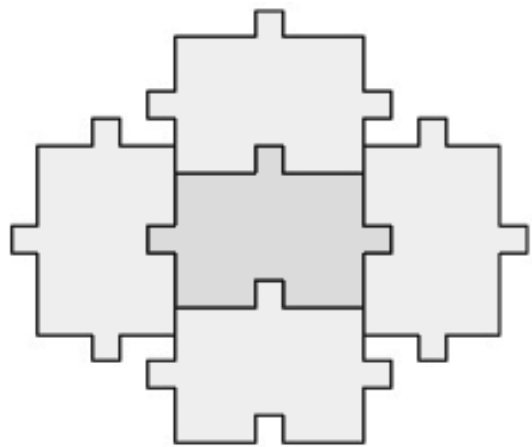

*Figure 4*

# A Region with Thickness Number = 4

**Lemma:** If a shape can form a wall of thickness 2n+1, it necessarily has Heesch number at least equal to n.

**Proof:** Indeed, a wall of thickness 2n+1 is, by definition, formed by 'stacking' without cavities 2n+1 walls of thickness 1 each. Consider any unit R in the $n^{th}$ of these thickness 1 walls. If we move from the boundary of R in any outward direction, we are guaranteed to pass through a minimum of n units before we exit the big wall of thickness 2n+1. This plus the requirement that the big wall has no cavities guarantees that within the wall, unit R sits surrounded by n layers of copies of itself. QED.

Searching for non-tiles which can form walls of thickness greater than 2, we find a non-convex figure – a hexagon with two small outward bulges and an inward dent that forms a wall of thickness 4 as shown in figure 2.

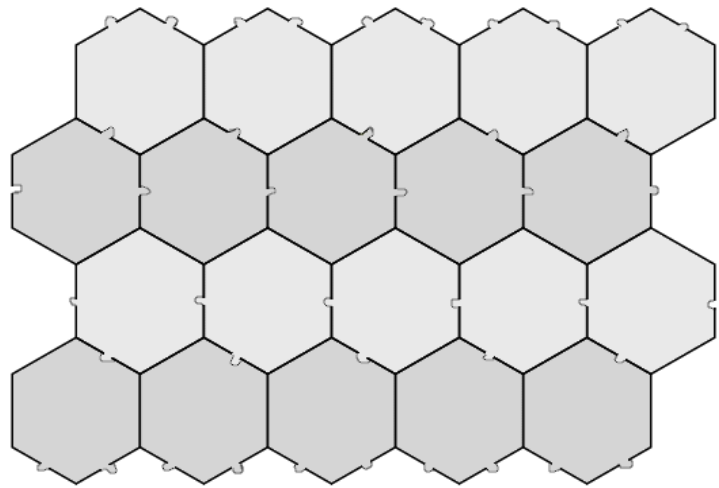

*Figure 5*

**Theorem:** The above deformed hexagon has thickness number 4.

**Proof:** It is sufficient to show that the hexagon cannot form walls of thickness 5 (that would automatically imply no higher thickness wall is possible). From the lemma above, for a wall of thickness 5 to be formed, the shape necessarily should have a Heesch number of at least 2 (ie. one should be able to surround a central unit with two layers of copies). So here, we only need to prove that the region cannot form two layers of copies of itself around a central unit.

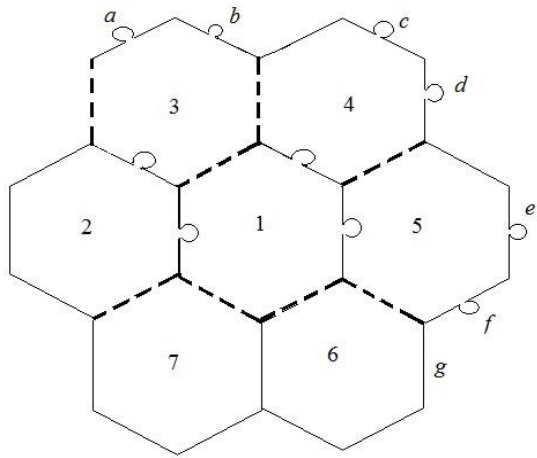

*Figure 6*

As in figure 6, let '1' be the central unit and 2-7 be the units surrounding it forming the first layer of a Heesch arrangement around 1. Since a unit has to have outward bulges on 2 adjacent edges, Unit 2 has to have an outward bulge into one of units 3 or 7. Without loss of generality, we take 2 to have a bulge into 3 - and the latter has to have a corresponding inward dent. Note that the edges in the layout shown as dashed thick lines in figure 6 have to be necessarily straight (indeed, the two edges of a unit adjacent to an edge with an inward dent have to be both straight). We have not shown the inward dent on unit 2 above because it is not important in the discussion.

Consider the outer boundary of first layer around unit 1. We observe that if at any concave vertex of this boundary, both boundary edges meeting there are non-straight with at least one having an outward bulge, we cannot fit a fresh unit from outside into that concave vertex and hence the formation of another layer around the layout becomes impossible. So, going counter-clockwise from unit 3 up to and including unit 5, we necessarily need the edges named *a,b,c,d,e* and *f* to have outward bulges.

Consider units 6 and 7. Since their upper boundaries are entirely composed of straight edges, we need to have their vertical edges to be all non-straight. So, the edge *g* is necessarily non-straight. Now, *f* and *g* are adjacent edges meeting at a concave vertex of the outer boundary of the layout with *f* having an outward bulge, so *g* too being non-straights implies no further outer layer is possible. The region fails to have a Heesch number 2 and thus fails to form walls of thickness 5 and above. QED

## Other Values of Thickness Number

We now show a non-convex shape (also discussed in [3]) that can form walls of thickness 3 (figure 7) with a conjecture that its thickness number is also 3.

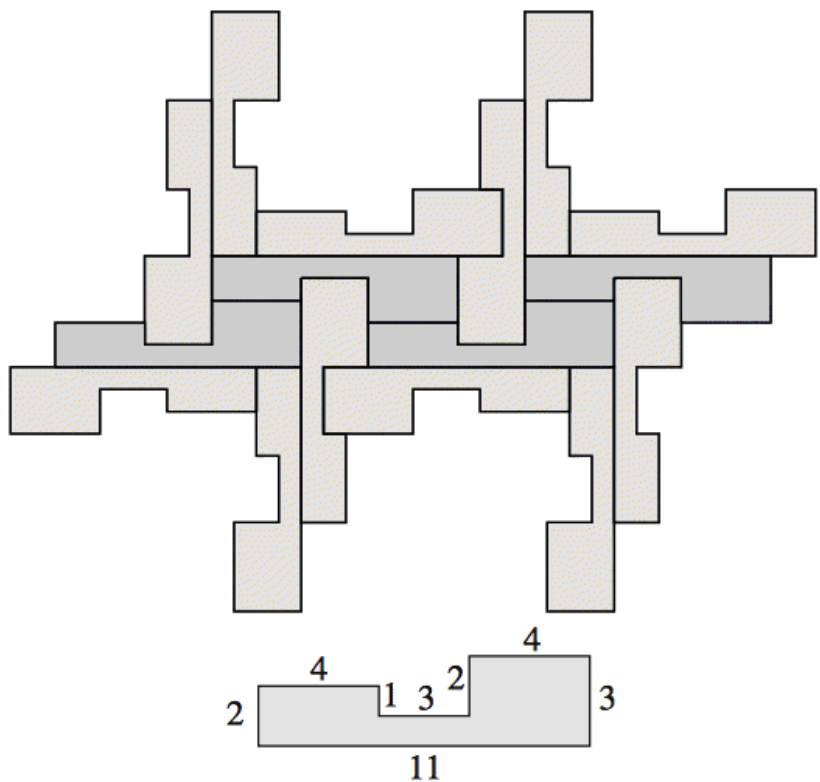

*Figure 7*

## Further Questions

1. Are there non-tile regions which can form walls of thickness greater than 4? If so, is there an upper limit on possible values of the thickness number for non-tiles?
2. Are there *convex* non-tile regions which can form walls of thickness greater than 2?
3. Are there convex regions with Heesch number 1 but for which thickness number is only 1?

# References

1. https://en.wikipedia.org/wiki/Heesch%27s_problem
2. Erich Friedman: Heesch Tiles with Surround Numbers 3 and 4, Geombinatorics, Volume VIII, 4 (1999) 101-103.
3. Casey Mann: Heesch's Tiling Problem, Amer. Math.. Monthly, Vol. Vol. 111 (6), June-July 2004, 509-517
4. https://www.uwgb.edu/dutchs/symmetry/heesch.htm